\documentclass[11pt]{article}
\usepackage{latexsym,amsmath}
\usepackage[psamsfonts]{amssymb}

\newtheorem{theorem}{Theorem}

\newtheorem{lemma}[theorem]{Lemma}
\newtheorem{proposition}[theorem]{Proposition}
\newtheorem{corollary}[theorem]{Corollary}

\def\endofproof{\nobreak\hfill $\blacksquare$ \goodbreak}

\def\Nmath{\mathbb{N}}
\def\N{\Nmath}

\def\Fmath{\mathbf{F}}

\def\dfn#1{\textbf{\emph{#1}}}  

\def\GG{\mathcal{G}}  
\def\EE{\mathtt{E}}  
\def\VV{\mathtt{V}}  
\def\origin{\rho} 

\def\RR{\mathcal{R}}
\def\SS{\mathcal{S}}
\def\cluster{\mathrm{cl}}
\def\RRcluster{\mathcal{R}^{\cluster}}
\def\UPsi{\overline\Psi}

\def\wsf{\mathcal{F}_1}
\def\fmsf{\mathcal{F}_2}

\def\cost{\mathcal{C}}

\def\Leb{\textit{Leb}}

\begin{document}
\thispagestyle{empty}
\title{A Measurable-Group-Theoretic Solution to von Neumann's Problem}
\author{
Damien Gaboriau\thanks{CNRS} \and 
Russell Lyons\thanks{Supported partially by NSF grant DMS-0705518 and Microsoft Research.}}
\date{\today}
\maketitle

\begin{abstract}
{We give a positive answer, in the measurable-group-theory context, to von
Neumann's problem of knowing whether a non-amenable countable discrete
group contains a non-cyclic free subgroup. We also get an embedding result
of the free-group von Neumann factor into restricted wreath product factors.}
\end{abstract}

\bigskip
\noindent
{\small \textbf{2000 Mathematical Subject Classification}: 37A20, 20E05, 82B43, 05C80}


\bigskip
\bigskip

Amenability of groups is a concept introduced by J.~von~Neumann in his seminal article \cite{vN29} to explain the so-called Banach-Tarski paradox. 
He proved that a discrete group containing the free group $\Fmath_{2}$ on
two generators as a subgroup is non-amenable. Knowing whether this was a
characterization of non-amenability became known as von Neumann's Problem
and was solved by the negative by A.~Ol$'${\v{s}}anski{\u\i} \cite{Ol80}.
Still, this characterization could become true after relaxing the notion of ``containing a subgroup''.
K.~Whyte gave a very satisfying geometric group-theoretic solution: \emph{
A finitely generated group $\Gamma$ is non-amenable iff it  admits a partition with pieces uniformly bilipschitz equivalent to the regular 4-valent tree  \cite{Wh99}.}
Geometric group theory admits a measurable counterpart, namely, measurable
group theory. The main goal of our note is to provide a solution to von Neumann's problem in this context. 
We show that any countable non-amenable group admits a measure-preserving free action on some probability space, such that the orbits may be measurably partitioned into pieces given by an $\Fmath_{2}$-action.

To be more precise, we use the following notation.
For a finite or
countable set $M$, let $\mu$  denote the product $\otimes_M \Leb$  on $[0,1]^{M}$ of the Lebesgue measures $\Leb$ on $[0,1]$, and for $p\in [0,1] $, let $\mu_p$ denote the
product of the discrete measures $(1-p)\delta_{\{0\}}+p\delta_{\{1\}}$ on
$\{0,1\}^{M}$.
Thus, the meaning of $\mu$ may vary from use to use as $M$ varies.
Usually $M$ will be a countable group $\Lambda$ or the set $\EE$ of edges of a Cayley graph.

\begin{theorem}\label{Th: main theorem}For any countable discrete  
non-amenable group $\Lambda$,
there is a measurable ergodic essentially free action of\/ $\Fmath_2$ on $([0,1]^ 
{\Lambda},\mu)$
such that almost every  $\Lambda$-orbit of the Bernoulli shift decomposes into $\Fmath_2$-orbits.
\end{theorem}
In other words, the orbit equivalence relation of the $\Fmath_{2}$-action is contained in that of the $\Lambda$-action. 
We give two proofs of this theorem, each with its own advantages.

For some purposes, it is useful to get a Bernoulli shift action with
a discrete base space. We show:
\begin{theorem}\label{Th: main theorem with discrete base space}
For any finitely generated non-amenable group $\Gamma$, there is $n\in  
\Nmath$ and a non-empty interval $(p_1,p_2)$ of parameters $p$ for  
which there is an ergodic essentially free action of $\Fmath_2$ on $\prod_ 
{1}^{n} (\{0,1\}^{\Gamma}, \mu_p)$
such that almost every  $\Gamma$-orbit of the diagonal Bernoulli shift decomposes into $\Fmath_2$-orbits.
\end{theorem}

These results have operator-algebra counterparts:
\begin{corollary}\label{Cor: wreath prod with infinite H}
Let $\Lambda$ be a countable discrete non-amenable group and $H$ be an  
infinite group. Then the von Neumann factor $L(H \wr \Lambda )$ of  
the restricted wreath product contains a copy of the von Neumann  
factor $L(\Fmath_2)$ of the free group.
\end{corollary}

\begin{corollary}\label{Cor: wreath prod with finite H}
Let $\Gamma$ be a finitely generated discrete non-amenable group.
Let $n$, $p_1$, $p_2$ be as in Theorem~\ref{Th: main theorem with  
discrete base space} and let $p=\frac{\alpha}{\beta}\in (p_1,p_2)$,  
with $\alpha, \beta\in \Nmath$.
Assume that $H$ contains an abelian subgroup $K$ of order $k=\beta^n 
$. Then the von Neumann factor $L(H \wr \Lambda )$ of the restricted  
wreath product contains a copy of the von Neumann factor $L(\Fmath_2) 
$ of the free group.
\end{corollary}

\medskip
For this paper, we assume a certain familiarity with the results and notation of 
 \cite{Gab05}, \cite{Gab00}  and \cite{LS99}.

\bigskip
\noindent
\textbf{{\large Acknowledgment}}\\
We are very grateful to Sorin Popa for bringing to our attention the above corollaries.
We also thank the referee for a careful reading.

\begin{centerline} {{------ }O{ ------}} \end{centerline}

A \dfn{(countable standard) equivalence relation} on the standard Borel
space $(X,\nu)$ is an equivalence relation $\RR$ with countable classes
that is a Borel subset of $X\times X$ for the product $\sigma$-algebra. 

 A \dfn{(measure-preserving oriented) graphing} on $(X,\nu)$ is a
denumerable family $\Phi=(\varphi_{i})_{i\in I}$ of partial
measure-preserving isomorphisms $\varphi_{i}:A_{i}\to B_{i}$ between Borel
subsets $A_{i},B_{i}\subset X$.

A graphing $\Phi$ \dfn{generates} an equivalence relation $\RR_{\Phi}$:
the smallest equivalence relation that contains all pairs $(x,
\varphi_i(x))$.
The \dfn{cost} of a graphing $\Phi=(\varphi_{i})_{i\in I}$
is the sum of the measures of the domains $\sum_{i\in I} \nu(A_{i})$.
The \dfn{cost}, $\mathrm{cost}(\RR, \nu)$, of 
$(\RR, \nu)$ is the infimum of the costs of the
graphings that generate $\RR$. 
The \dfn{graph (structure)} $\Phi[x]$ of a graphing $\Phi$ at a point $x
\in X$ is the graph whose vertex set is the equivalence class $\RR_{\Phi}[x]$ of $x$
and whose edges are the pairs $(y, z) \in \RR_{\Phi}[x] \times \RR_{\Phi}[x]$ such that for
some $i \in I$, either $\varphi_i(y) = z$ or $\varphi_i(z) = y$.
For more on cost, see \cite{Gab00} or \cite{Kechris-Miller}.

\begin{centerline} {{------ }O{ ------}} \end{centerline}

\bigskip
\noindent
\textbf{{\large Proofs}}

Since the union of an increasing sequence of amenable groups is still  
amenable,
$\Lambda$ contains a non-amenable finitely generated subgroup. Let $ 
\Gamma$ be such a subgroup.

\begin{centerline} {{------ }O{ ------}} \end{centerline}

If $S$ is a finite generating set of $\Gamma$ (maybe with  
repetitions), $\GG=(\VV,\EE)$ denotes the associated right Cayley graph (with  
vertex set $\VV$): The set $\EE$ of edges is indexed by $S$ and $\Gamma$. Given $s
\in S$ and $\gamma \in \Gamma$, the corresponding edge
is oriented from the vertex $\gamma$ to $\gamma s$. 
Note that $\Gamma$ acts freely on $\GG$ by multiplication on the left. Let $\origin:=\mathrm{id}$, the identity of the group $\Gamma$, chosen as base vertex for $\GG$.

The set of the subgraphs of $\GG$ (with the same set of  
vertices $\VV$) is naturally identified with $\Omega:=\{0,1\}^{\EE}$.  
The connected components of $\omega\in \Omega$ are called its \dfn
{clusters}.

Consider a probability-measure-preserving essentially free (left) $\Gamma$-action
on some standard Borel space $(X,\nu)$ together with a $\Gamma$-equivariant
Borel map $\pi:X\to \{0,1\}^{\EE}$.

The \dfn{full} equivalence relation $\RR_{\Gamma}$ generated by the
$\Gamma$-action $X$ is graphed by the graphing $\Phi=(\varphi_{s})_{s\in
S}$,  where  $\varphi_{s}$ denotes the action by $s^{-1}$.

We define the following equivalence subrelation on $X$ (see
\cite[Sect.~1]{Gab05}): the \dfn 
{cluster equivalence subrelation} $\RRcluster$, graphed by the  
graphing $\Phi^{\cluster}:=(\varphi^{\cluster}_{s})_{s\in S}$ of  
partial isomorphisms, where $\varphi^{\cluster}_{s}$ is 
the restriction $\varphi^{\cluster} 
_s:=\varphi_s\vert A_s$ of $\varphi_s$ to the Borel subset $A_s$ of $x\in X$  
for which the edge $e$ labelled $s$ from $\origin$ to $\origin s$ 
lies in $\pi(x)$, i.e., $\pi(x)(e)=1$. 
Consequently, \emph{two points $x,y\in X$ are $\RRcluster$-equivalent
if and only if there is some $\gamma\in \Gamma$ such that $ 
\gamma^{-1} x= y$ and the vertices $\origin, \gamma\origin $ are in  
the same cluster of $\pi(x)$.}

The graph structure $\Phi^{\cluster}[x]$ given by the graphing $\Phi^ 
{\cluster}$ to the $\RRcluster$-class of any $x\in X$ is naturally  
isomorphic with the cluster $\pi(x)_\origin$ of $\pi(x)$ that  
contains the base vertex.
Denote by $U^{\infty}\subset X$ the Borel set of $x\in X$ whose $ 
\RRcluster$-class is infinite
and by $\RR_{\Gamma\vert \infty}$ (resp. $\RRcluster_{\vert\infty}$) the restriction of $\RR_{\Gamma}$ (resp. $\RRcluster$) to $U^{\infty}$.

Write $\mathcal{P}(Y)$ for the power set of $Y$.
The map $X\times \VV \to X$ defined by $(x,\gamma\origin)\mapsto \gamma^{-1}.x$ induces a map 
$\Psi : X\times \mathcal{P}(\VV)\to \mathcal{P}(X)$ 
that is invariant under the (left) diagonal $\Gamma$-action (i.e.,
$\Psi(\gamma. x, \gamma. C) = \Psi(x, C)$ for all $\gamma \in \Gamma$, $x \in
X$, and $C \subset \VV$)
and such that $\Psi(x,\VV)$ is the whole $\RR_{\Gamma}$-class of $x$.
The restriction of $\Psi$ to the $\Gamma$-invariant subset $\mathfrak{C}_{\infty}^{\cluster}:=\bigl\{(x,C): x\in X, C\in \mathcal{P}(\VV), C \mathrm{\ is\ an\ infinite\ cluster\ of\ }\pi(x)
\bigr\}$
 sends each $(x,C)$ (and its $\Gamma$-orbit) to a whole infinite
$\RRcluster$-class, namely, the $\RRcluster$-class of $\gamma^{-1}. x$ for any $\gamma$ such that  $\gamma\origin\in C$. Moreover, for each 
$x\in U^{\infty}$, its $\RR_{\Gamma\vert\infty}$-class decomposes into infinite $\RRcluster$-sub-classes that are in one-to-one correspondence with the elements of 
$ \mathfrak{C}_{\infty}^{\cluster}$ that have $x$ as first coordinate.
Note that the set $\{(x, y, C) \in X \times X \times \mathcal P(\VV) :
x \in \Psi(y, C)\}$ is Borel, whence for a
Borel set $\mathcal A\subset \mathfrak{C}_{\infty}^{\cluster}$, the set
$\UPsi(\mathcal A) := \bigcup \Psi(\mathcal A)$ is measurable, being the
projection onto the first coordinate of the Borel set 
$\{(x, y, C) : x \in \Psi(y, C)\} \cap (X \times \mathcal A)$.

We say that $(\nu, \pi)$ has \dfn{indistinguishable infinite clusters}
if for every $\Gamma$-invariant Borel subset $\mathcal{A}\subset
\mathfrak{C}_{\infty}^{\cluster}$, the set of $x\in X$ for which some
$(x,C)\in \mathcal{A}$ and some $(x,C)\in \mathfrak{C}_{\infty}^{\cluster}
\setminus \mathcal{A}$ has $\nu$-measure $0$.
In other words, the $\RRcluster$-invariant partition
$U^{\infty}=\UPsi(\mathcal{A})\cup \complement\UPsi(\mathcal{A})$ is not allowed to split any 
$\RR_{\Gamma\vert \infty}$-class (up to a union of measure 0 of such classes).
The following proposition, using this refined notion of indistinguishability, corrects \cite[Rem.~2.3]{Gab05}.

\begin{proposition}
\label{Prop: ergodicity}
Let $\Gamma$ act ergodically on $(X,\nu)$ and 
$\pi:X\to \{0,1\}^{\EE}$ be a $\Gamma$-equivariant Borel map
such that $\nu(U^\infty) \ne 0$.
Then $\RRcluster_{\vert\infty}$, the cluster equivalence relation
restricted to its infinite locus $U^{\infty}$, is ergodic if and only if
$(\nu, \pi)$ has indistinguishable infinite clusters.
\end{proposition}

\noindent{\it Proof.}\enspace
Suppose that $\RRcluster_{\vert\infty}$ is ergodic.
Then for every $\Gamma$-invariant Borel subset $\mathcal{A}\subset
\mathfrak{C}_{\infty}^{\cluster}$, its image $\UPsi(\mathcal{A})$ is a union
of $\RRcluster_{\vert\infty}$-classes, whence $\UPsi(\mathcal{A})$ or its
complement $\complement\UPsi(\mathcal{A})$ in $U^\infty$ has measure 0.
In particular, the partition $U^{\infty}=\UPsi(\mathcal{A})\cup
\complement\UPsi(\mathcal{A})$ is trivial, whence
$\nu$ has indistinguishable infinite clusters.

Conversely, suppose that
$(\nu, \pi)$ has indistinguishable infinite clusters.
An $\RRcluster_{\vert\infty}$-invariant partition $U^{\infty} =
\mathcal{U}\cup \complement\mathcal{U}$ defines a partition
$\mathfrak{C}_{\infty}^{\cluster} = \mathcal{A}\cup \complement\mathcal{A}$
according to whether $\Psi(x,C)\in \mathcal{U}$ or
$\complement\mathcal{U}$.
Then for $\nu$-almost every $x\in U^\infty$, all $\Psi(x,C)$ are in
$\mathcal{U}$ or all are in its complement, i.e., the
$\RRcluster_{\vert\infty}$-subclasses into which the
$\RR_{\Gamma\vert\infty}$-class of $x$ splits all belong to one side.
Since $\RR_{\Gamma\vert\infty}$ is $\nu$-ergodic, this side has to be the
same for almost every $x$.
This means that the other side is a null set.
This holds for any partition $\mathcal{U}\cup \complement\mathcal{U}$,
whence
$\RRcluster_{\vert\infty}$ is ergodic. 
\endofproof

\bigskip

If $X$ has the form $X = \Omega \times Y$, then $\nu$ is called
\dfn{insertion tolerant} (see \cite{LS99}) if for each edge $e\in \EE$,
the map $\Pi_e:X\to X$ defined by $(\omega, y) \mapsto \big(\omega\cup
\{e\}, y\big)$ quasi-preserves the measure, i.e., $\nu(A)>0$
implies $\nu\big(\Pi_{e}(A)\big)>0$ for every measurable  subset
$A\subseteq X$.
Call a map $\pi : X \to \Omega$ \dfn{increasing} if
$\pi(\omega, y) \supseteq \omega$ for all $\omega \in \Omega$.
An action of $\Gamma$ on $\Omega \times Y$ is always assumed to act on the
first coordinate in the usual way.
A slight extension of \cite[Th.~3.3, Rem.~3.4]{LS99}, proved in the same
way, is the following:

\begin{proposition}
\label{Prop: indistinguishability}
Assume that $\Gamma$ acts on $(\Omega \times Y,\nu)$ preserving the measure
and $\pi:\Omega \times Y\to \Omega$ is an increasing $\Gamma$-equivariant
Borel map with $\nu(U^\infty)\ne 0$.
If $\nu$ is insertion tolerant, then $(\nu, \pi)$ has indistinguishable
infinite clusters.
\end{proposition}

\begin{proposition}
\label{Prop: GamLam}
If\/ $\Gamma < \Lambda$, then there are $\Gamma$-equivariant isomorphisms $([0,1]^{\EE},\mu)\simeq ([0,1]^\Gamma,\mu)\simeq ([0,1]^{\Lambda},\mu)$
between the Bernoulli shift actions of\/ $\Gamma$. In particular, the orbits of the Bernoulli shift $\Lambda$-action on $[0,1]^{\Lambda}$ are partitioned into subsets that are identified with the orbits of the standard Bernoulli shift $\Gamma$-action on $[0,1]^{\Gamma}$.
\end{proposition}

\noindent{\it Proof.}\enspace
A countable set $\EE$ on which $\Gamma$ acts freely may be decomposed  
by choosing a representative in each orbit so as to be identified with  
a disjoint union of $\Gamma$-copies, $\EE\simeq \coprod_{J} \Gamma$,
and to give $\Gamma$-equivariant identifications
$[0,1]^{\EE}=[0,1]^{\coprod_{J} \Gamma}=([0,1]^{J})^{\Gamma}$.

The edge set $\EE\simeq \coprod_{S} \Gamma$ of the Cayley graph of $ 
\Gamma$, as well as $\Lambda\simeq \coprod_{I} \Gamma$ 
 are such countable $\Gamma$-sets. 
Then  
isomorphisms of standard Borel probability spaces $([0,1], \Leb)
\simeq ([0,1]^{S}, \otimes_S\Leb)\simeq ([0,1]^{I},\otimes_I \Leb)$ induce $ 
\Gamma$-equivariant isomorphisms of the Bernoulli shifts:
\[
\begin{array}{ccccc}
[0,1]^{\Gamma} & \simeq & ([0,1]^{S})^{\Gamma} &\simeq &([0,1]^{I})^ 
{\Gamma} \\
\parallel & & \parallel & & \parallel \\
   {[0,1]^{\Gamma}} & \simeq & [0,1]^{\EE} &\simeq & [0,1]^{\Lambda}
\,.
\end{array}
\]
\endofproof

\begin{centerline} {{------ }O{ ------}} \end{centerline}

A subgraph $(\VV', \EE')$ of a graph $(\VV, \EE)$ is called \dfn{spanning}
if $\VV' = \VV$.
A vertex $a$ in a graph is called a \dfn{cutvertex} if there are two other
vertices in its component with the property that every path joining them
passes through $a$.
A \dfn{block} of a graph is a maximal connected subgraph that has no
cutvertex.
Every simple cycle of a graph is contained within one of its blocks.

\begin{lemma}
\label{Lem: blocks}
If all vertices of a block have finite degree and for each pair of vertices
$(a, b)$ there are only finitely many distinct paths joining $a$ to $b$,
then the block is finite.
\end{lemma}

\noindent{\it Proof.}\enspace
Suppose for a contradiction that the block is infinite.
Then it contains a simple infinite path $P$ of vertices $a_1, a_2, \dots$.
By Menger's theorem, $a_1$ and $a_n$ belong to a simple cycle $C_n$
for each $n > 1$.
But this implies that there are infinitely many distinct paths joining
$a_1$ to $a_2$: 
Fix $n$ and
let $a_j$ ($2 \le j \le n$) be the vertex on $P \cap C_n$ with minimal
index $j$.
We may assume that $C_n$ is oriented so that it visits
$a_n$ before it visits $a_j$.
Then simply follow $C_n$ from $a_1$ until $a_j$,
and then follow $P$ to $a_2$.
\endofproof

\begin{proposition}(For any Cayley graph)
\label{Prop: forest}
Let $X := \Omega \times [0, 1]^\Gamma$ and $\epsilon > 0$.
Let $\nu := \mu_\epsilon \times \mu$.
There is a $\Gamma$-equivariant Borel map $f: X\to \Omega$
such that 
$(\nu, f)$ has indistinguishable infinite clusters
and
for all sufficiently small $\epsilon$,
the push-forward measure $f_{*}(\nu)$ of $\nu$ is  
supported on the set of spanning subgraphs of $\GG$
each of whose components is a tree with infinitely many ends.
\end{proposition}

\noindent{\it Proof.}\enspace
We may equivariantly identify $\big([0, 1]^\Gamma, \mu\big)$ with $\big([0,
1]^{\N \times \Gamma} \times [0, 1]^\EE, \mu \times \mu\big)$, so we
identify $(X, \nu)$ with
$\big(\Omega \times [0, 1]^{\N \times \Gamma} \times [0, 1]^\EE,
\mu_\epsilon \times \mu \times \mu\big)$.
Fix an ordering of $S \coprod S^{-1}$; this determines an
ordering of the edges incident to each vertex in $\GG$, where we ignore
edge orientations for the rest of this proof.
With $d$ denoting the degree of $\GG$, define the function $h(t) :=
\lceil{d t}\rceil$ for $t \in [0, 1]$.
Given a point $x = (\omega, ({r}(n, \gamma))_{n \in \N, \gamma \in \Gamma},
u) \in X$, construct the wired spanning forest $\wsf$ of $\GG$
by using the cycle-popping algorithm of D.~Wilson \cite[Sect.~3]{Wilson1996} as adapted in \cite[Th.~5.1]{BLPS01}, also called there ``Wilson's algorithm rooted at infinity'',
applied to the stacks where the $n$th edge in the stack under $\gamma$ is
defined as the $h({r}(n, \gamma))$th edge incident to $\gamma$.
The measure $\nu$ is insertion tolerant and the
map $\pi: x \mapsto \omega \cup \wsf$ is increasing, whence by Proposition
\ref{Prop: indistinguishability}, the pair $(\nu, \pi)$ has
indistinguishable infinite clusters.
Notice that all clusters are infinite.
Now use $u$ to construct the free minimal spanning forest $\fmsf$ in each
cluster of $\pi(x)$, that is, for every cycle $\Delta \subset \pi(x)$,
delete the edge $e \in \Delta$ with maximum $u(e)$ in that cycle.
The map $f$ is $f(x) := \fmsf$.

Now the $\nu$-expected number
of distinct simple paths in $\pi(x)$
that join any two vertices is finite (equation
(13.7) of \cite{BLPS01}) for all sufficiently small $\epsilon$. 
In particular, the number of such paths is finite $\nu$-a.s. 
By Lemma \ref{Lem: blocks},
this means that all blocks of $\pi(x)$ are finite, so that $\fmsf$ is a
spanning tree in each block.
Therefore each component of $\fmsf$ spans a component of $\pi(x)$.
Thus, $\fmsf$ determines the same cluster relation and so $(\nu, f)$ also has
indistinguishable (infinite) clusters.
Finally, the fact that the clusters of $\pi(x)$, and hence those of
    $\fmsf$, have infinitely many ends follows, e.g., from
    \cite[Th.~13.7]{BLPS01}.
\endofproof

\medskip
The cluster relation determined by $f$ of Proposition \ref{Prop:
forest} is treeable and has cost larger than 1 by \cite[Cor. IV.24
(2)]{Gab00}, has finite cost (since the degree is bounded), and is ergodic by Proposition \ref{Prop:
ergodicity}. 
Since we may equivariantly identify $\big(\Omega \times [0, 1]^\Gamma,
\mu_\epsilon \times \mu\big)$ with $\big([0, 1]^\Gamma, \mu\big)$, 
we proved:
\begin{proposition}
\label{Prop: BLPSforest}
For any Cayley graph of $\Gamma$, the Bernoulli action on
$([0,1]^{\Gamma},\mu)$ contains a treeable
subrelation that is ergodic and has cost in the interval $(1,\infty)$.
\end{proposition}

At this point, we already have a reasonable answer to the analogue  
of von Neumann's problem, since ``treeable relation'' is the analogue of  
``free group'' and cost $\cost>1$ is, in the context of treeable relations,
equivalent to non-amenability.

\begin{centerline} {{------ }O{ ------}} \end{centerline}

\medskip
An alternative approach begins with a more
explicit $f$ and a more common  
measure $f_{*}(\mu)$, namely, the Bernoulli measure $\mu_p$ on $\{0,1\} 
^{\EE}$ for a certain parameter $p$, but requires us to choose a  
particular Cayley graph for $\Gamma$.
It also requires us to obtain a treeable subrelation in a less explicit
way.
This is accomplished as follows.

Results of H{\"a}ggstr{\"o}m-Peres \cite{HP99} imply that there are two
critical values $0 < p_c \le p_u \le 1$ such that
\\
-  (finite phase, $p\in [0,p_c)$) $\mu_p$-a.s., the subgraph has only finite  
clusters; \\
- (\dfn{non-uniqueness phase}, $p\in (p_c,p_u)$) $\mu_p$-a.s., infinitely
  many of the clusters of the subgraph are
infinite, each one with  
infinitely many ends; \\
-  (the uniqueness phase, $p\in (p_u,1]$) $\mu_p$-a.s., the subgraph has only one  
cluster that is infinite. 

The situation for the critical values $p_c$ and $p_u$ themselves is  
far from clear. Benjamini and Schramm \cite{BS96} conjectured that  
$p_c\not= p_u$ for every Cayley graph of a f.g.\ non-amenable group.  The
main result of \cite{PS00} (Th.~1, p.~498) asserts that given a f.g.\
non-amenable group $ 
\Gamma$, there is a finite set of generators such that the associated  
Cayley graph admits a non-trivial interval of non-uniqueness.
Thus:
\begin{proposition}(For particular Cayley graphs)
\label{Prop: Bernoulli perc. non uniqueness}
There exists a Cayley graph of $\Gamma$ and  
a non-empty interval $(p_c,p_u)$ such that, for any $p\in (p_c,p_u)$,  
the Bernoulli measure $\mu_p$ on $\{0,1\}^{\EE}$ is  supported on the  
set of subgraphs admitting infinite components, each one with  
infinitely many ends.
\end{proposition}

Let $\pi:(X, \nu)\to \{0,1\}^{\EE}$ denote
either\\
   (i) $f_p: ([0,1]^{\EE},\mu)\to \{0,1\}^{\EE}$ induced by the characteristic function $\chi_{[0,p]}: [0,1]\to \{0,1\}$ of 
   $[0,p]$,
or \\
   (ii) the identity map $(\{0,1\}^{\EE}, \mu_p)\to \{0,1\}^{\EE}$,\\
both with the natural Bernoulli $\Gamma$-action. Notice that the action is essentially free when $0 < p <1$.

In case (ii), we have that $(\mu_p, \pi)$ has indistinguishable infinite
clusters by \cite[Th.~3.3]{LS99}.
Case (i) is essentially the same, but first we must
identify $([0,1]^{\EE},\mu)$ equivariantly as $(\{0,1\}\times
[0,1])^{\EE}=\{0,1\}^{\EE}\times [0,1]^{\EE}$ equipped with the product measure $\mu_{p}\times \mu$ in such a way that $f_p$ becomes the identity on the first coordinate.
Then we have insertion tolerance and so, by \cite[Rem.~3.4]{LS99},
indistinguishable infinite clusters.

Hence, in both cases,
for any $p$ given by Prop.~\ref{Prop: Bernoulli perc. non uniqueness}, the
locus $U^{\infty}$ of infinite classes of $\RRcluster$ is non-null
and we have
ergodicity of the restriction $\RRcluster_{\vert\infty}$ of $\RRcluster$ to $U^{\infty}$ by Proposition~\ref{Prop: ergodicity}. We claim that its normalized cost (i.e., computed with respect to the normalized probability measure $\nu/\nu(U^{\infty})$ on $U^{\infty}$) satisfies $1<\cost(\RRcluster_{\vert\infty})<\infty$. The finiteness of the cost is clear since $S$, the index set for $\Phi^{\cluster}$, is finite. That it is strictly greater than $1$ is a direct application of \cite[Cor. IV.24 (2)]{Gab00}, since the graph $\Phi^{\cluster}[x]\simeq \pi(x)_\origin$ associated with almost every $x\in U^{\infty}$ has at least $3$ ends.

In order to extend $\RRcluster_{\vert\infty}$ to a
subrelation of $\RR_{\Gamma}$ defined on the whole of $X$,
choose an enumeration $\{\gamma_i\}_{i\in \Nmath}$ of $\Gamma$. For each $x\in X\setminus U^{\infty}$, let $\gamma_x$ be the first element $\gamma_j\in \Gamma$ such that $\gamma_j\cdot x\in U^{\infty}$. 
Then the smallest equivalence relation containing
$\RRcluster_{\vert\infty}$ and the $(x,\gamma_x\cdot x)$'s is a subrelation of $\RR_\Gamma$, is ergodic, and has cost in $(1, \infty)$
by the induction formula of \cite[Prop.
II.6]{Gab00}. We proved:

\begin{proposition}\label{prop:subrel erg+cost}
For a Cayley graph and a $p$ given by Proposition~\ref{Prop: Bernoulli
perc. non uniqueness}, the Bernoulli actions on both $([0,1]^{\EE},\mu)$ and $(\{0,1\}^{\EE},\mu_{p})$ contain a subrelation that is ergodic and has cost in the open interval $(1,\infty)$.
\end{proposition}

\begin{centerline} {{------ }O{ ------}} \end{centerline}

\begin{proposition}\label{prop: contains treeable}
If an equivalence relation $\RR$
is ergodic and has cost in $(1,\infty)$, then it contains a treeable subrelation $\SS_{1}$ that is ergodic and has cost in $(1,\infty)$.
\end{proposition}
\noindent{\it Proof.}\enspace
This is ensured by a result proved independently by A.~Kechris and
B.~Miller~\cite[Lem.~28.11; 28.12]{Kechris-Miller} and by M.~Pichot~\cite 
[Cor.~40]{Pichot-Thesis}, through a process of erasing cycles from a  
graphing of $\SS_1$ with finite cost that contains an ergodic global  
isomorphism.\endofproof

\begin{centerline} {{------ }O{ ------}} \end{centerline}
\begin{proposition}\label{prop: erg F2 action}
If a treeable equivalence relation $\SS_{1}$
is ergodic and has cost $> 1$, then it contains a.e.\ a subrelation $\SS_{2}$ that is generated by an ergodic free action of the free group $\Fmath_2$.
\end{proposition}
\noindent{\it Proof.}\enspace
If the cost of $\SS_1$ is $>2$, this follows from a result of  
G.~Hjorth~\cite{Hjorth06} (see also \cite[Sect.~28]{Kechris-Miller}).
Otherwise, one first considers the restriction of the treeable $\SS_1
$ to a small enough Borel subset $V$: this increases the normalized  
cost by the induction formula of \cite[Prop. II.6 (2)]{Gab00} to get  
$\cost(\SS_1\vert V)\geq 2$.
In fact, it follows from the proof of \cite[Th.~28.3]{Kechris-Miller}  
that
\emph{a treeable probability measure-preserving equivalence relation  
with cost $\geq 2$ contains a.e.\ an equivalence subrelation that is given  
by a free action of the free group $\Fmath_2=\langle a,b\rangle$ in  
such a way that the generator $a$ acts ergodically}. By considering a  
subgroup of $\Fmath_2$ generated by $a$ and some conjugates of $a$ of  
the form $b^{k}ab^{-k}$, one gets an ergodic treeable subrelation of $ 
\SS_1\vert V$ with a big enough normalized cost that, when extended  
to the whole of $X$ (by using partial isomorphisms of $\SS_1$),
it gets cost $\geq 2$ (by the induction back \cite[Prop. II.6 (2)] 
{Gab00}) and of course remains ergodic.
Another application of the above-italicized result gives the desired  
ergodic action of $\Fmath_2$ on $X$.
\endofproof

\begin{centerline} {{------ }O{ ------}} \end{centerline}

The proof of Theorem~\ref{Th: main theorem with discrete base space} is now
complete as a direct consequence of Propositions~\ref{prop:subrel erg+cost}
(for the case $X = \{0, 1\}^\EE$), \ref{prop: contains treeable} and \ref{prop: erg F2 action}.
\endofproof

\begin{centerline} {{------ }O{ ------}} \end{centerline}

In case $X=[0,1]^{\EE}$ of Prop. \ref{prop:subrel erg+cost}, by using
Prop.~\ref{Prop: GamLam},
we can see $\SS_2$ (with $\SS_{2}\subset \SS_{1}\subset \RR_{\Gamma}$ given by Prop.~\ref{prop: erg F2 action} and \ref{prop: contains treeable}) as an equivalence subrelation of that given by the  
Bernoulli shift action of $\Lambda$. This finishes the proof of
Theorem~\ref{Th: main theorem}.
Alternatively, we may use Prop.~\ref{Prop: BLPSforest}
and a similar argument
to prove Theorem~\ref{Th: main theorem}.
\endofproof

\begin{centerline} {{------ }O{ ------}} \end{centerline}

\noindent{\it Proof of Cor~\ref{Cor: wreath prod with infinite H}.}\enspace
For any diffuse abelian subalgebra $A$ of $L(H)$, the von Neumann  
factor $L(H \wr \Lambda )=L(\Lambda \ltimes \oplus_{\Lambda} H)$  
contains the von Neumann algebra crossed product 
$\Lambda\ltimes\otimes_{\Lambda} A$,
which is isomorphic with the group-measure-space factor $\Lambda 
\ltimes L^{\infty}([0,1]^{\Lambda}, \mu)$ associated with the  
Bernoulli shift. The corollary then follows from Th.~\ref{Th: main  
theorem}.\endofproof

\begin{centerline} {{------ }O{ ------}} \end{centerline}

\noindent{\it Proof of Cor.~\ref{Cor: wreath prod with finite H}.}\enspace
If $\hat{K}$ is the dual group of $K$, then $L(H \wr \Gamma )$  
contains $L(K \wr \Gamma )$, which is isomorphic with the
group-measure-space factor $\Gamma\ltimes L^{\infty}(\hat{K}^{\Gamma})$
associated with the Bernoulli shift of $\Gamma$ on $\hat{K}^{\Gamma} 
$, where the finite set $\hat{K}\simeq \{1,2, \ldots, k\}$ is equipped  
with the equiprobability measure $\nu$.
The result is then obtained by taking the pull-back of the
$\Fmath_2$-action on $\prod_1^n \{0,1\}^{\Gamma}$, given in Th.~\ref{Th: main  
theorem with discrete base space}, by the $\Gamma$-equivariant Borel  
map $\hat{K}^{\Gamma}\to (\{0,1\}^{n})^{\Gamma}\simeq \prod_1^n \{0,1 
\}^{\Gamma}$, sending $\otimes \nu$ to $\mu_p$, that extends a map
$\{1,2,\ldots, k\}\to \{0,1\}^n$ (whose existence is ensured by the  
form of $k=\beta^n$).
\endofproof

\begin{centerline} {{------ }O{ ------}} \end{centerline}

It is likely that the free minimal spanning forest (FMSF) of a Cayley graph
of $\Gamma$ would serve as the desired ergodic subrelation $\SS_1$ of
Prop.~\ref{prop: contains treeable}, but its
indistinguishability, conjectured in \cite{LPS:msf}, is not known.
Also, it is not known to have cost $> 1$, but this is equivalent to $p_c <
p_u$, which is conjectured to hold and which we know holds for some Cayley
graph.
See \cite{LPS:msf} for information on
the FMSF and \cite{Tim06} for a weak form of indistinguishability.

\medskip
A general question remains open:\\
\textbf{Question:} Does every probability-measure-preserving free ergodic
action of a non-amenable countable group contain an ergodic subrelation generated by a free action of a non-cyclic free group?
More generally: Does every standard countable probability-measure-preserving
non-amenable ergodic equivalence relation contain a treeable non-amenable
ergodic equivalence subrelation?


\def\cprime{$'$}

\bigskip
{\footnotesize \hskip-\parindent Damien Gaboriau: \\
Unit\'e de Math\'ematiques Pures et Appliqu\'ees\\
Universit\'e de Lyon, CNRS, Ens-Lyon\\
46, all\'ee d'Italie \\
69364 Lyon cedex 7, France\\
{gaboriau@umpa.ens-lyon.fr}}

\bigskip

{\footnotesize \hskip-\parindent Russell Lyons: \\
 Department of Mathematics\\
Indiana University
Bloomington, IN 47405-5701 \\
USA\\
{rdlyons@indiana.edu}}

\end{document}